\documentclass[preprint,amsmath,amssymb,aps,showpacs]{revtex4-1} 
\usepackage{graphicx}
\usepackage{amssymb}
\usepackage{amsmath}
\usepackage{amsbsy}
\usepackage{bm}
\usepackage{tikz} 

\begin{document}
\setcounter{page}{1}
\newtheorem{t1}{Theorem}[section]
\newtheorem{d1}{Definition}[section]
\newtheorem{c1}{Corollary}[section]
\newtheorem{l1}{Lemma}[section]
\newtheorem{r1}{Remark}[section]

\newcommand{\cA}{{\cal A}}
\newcommand{\cB}{{\cal B}}
\newcommand{\cC}{{\cal C}}
\newcommand{\cD}{{\cal D}}
\newcommand{\cE}{{\cal E}}
\newcommand{\cF}{{\cal F}}
\newcommand{\cG}{{\cal G}}
\newcommand{\cH}{{\cal H}}
\newcommand{\cI}{{\cal I}}
\newcommand{\cJ}{{\cal J}}
\newcommand{\cK}{{\cal K}}
\newcommand{\cL}{{\cal L}}
\newcommand{\cM}{{\cal M}}
\newcommand{\cN}{{\cal N}}
\newcommand{\cO}{{\cal O}}
\newcommand{\cP}{{\cal P}}
\newcommand{\cQ}{{\cal Q}}
\newcommand{\cR}{{\cal R}}
\newcommand{\cS}{{\cal S}}
\newcommand{\cT}{{\cal T}}
\newcommand{\cU}{{\cal U}}
\newcommand{\cV}{{\cal V}}
\newcommand{\cX}{{\cal X}}
\newcommand{\cW}{{\cal W}}
\newcommand{\cY}{{\cal Y}}
\newcommand{\cZ}{{\cal Z}}

\def\cl{\centerline}
\def\bd{\begin{description}}
\def\be{\begin{enumerate}}
\def\ben{\begin{equation}}

\def\een{\end{equation}}

\def\benr{\begin{eqnarray}}
\def\eenr{\end{eqnarray}}
\def\benrr{\begin{eqnarray*}}
\def\eenrr{\end{eqnarray*}}
\def\bes {\begin{equation*}}
\def\ees {\end{equation*}}
\def\bbm {\begin{bmatrix}}
\def\ebm {\end{bmatrix}}
\def\ed{\end{description}}
\def\ee{\end{enumerate}}
\def\al{\alpha}
\def\b{\beta}
\def\bR{\bar\R}
\def\bc{\begin{center}}
\def\ec{\end{center}}
\def\d{\dot}
\def\D{\Delta}
\def\del{\delta}
\def\ep{\epsilon}
\def\g{\gamma}
\def\G{\Gamma}
\def\h{\hat}
\def\iny{\infty}
\def\La{\Longrightarrow}
\def\la{\lambda}
\def\m{\mu}
\def\n{\nu}
\def\noi{\noindent}
\def\Om{\Omega}
\def\om{\omega}
\def\p{\psi}
\def\pr{\prime}
\def\r{\ref}
\def\R{{\bf R}}
\def\ra{\rightarrow}
\def\s{\sum_{i=1}^n}
\def\si{\sigma}
\def\Si{\Sigma}
\def\t{\tau}
\def\th{\theta}
\def\Th{\Theta}

\def\vep{\varepsilon}
\def\vp{\varphi}
\def\pa{\partial}
\def\un{\underline}
\def\ov{\overline}
\def\fr{\frac}
\def\sq{\sqrt}

\def\WW{\begin{stack}{\circle \\ W}\end{stack}}
\def\ww{\begin{stack}{\circle \\ w}\end{stack}}
\def\st{\stackrel}
\def\Ra{\Rightarrow}
\def\R{{\mathbb R}}
\def\bi{\begin{itemize}}
\def\ei{\end{itemize}}
\def\i{\item}
\def\bt{\begin{tabular}}
\def\et{\end{tabular}}
\def\lf{\leftarrow}
\def\nn{\nonumber}
\def\va{\vartheta}
\def\wh{\widehat}
\def\vs{\vspace}
\def\Lam{\Lambda}
\def\sm{\setminus}
\def\ba{\begin{array}}
\def\ea{\end{array}}
\def\bd{\begin{description}}
\def\ed{\end{description}}
\def\lan{\langle}
\def\ran{\rangle}
\def\mf{\mathbf}
\def\ts{\times}
\def\ots{\otimes}
\def\bs{\boldsymbol}
\def\l{\label}
\def\r{\ref}
\numberwithin{equation}{section}

\title{On Gracefully Labeling Trees}

\author{Dhananjay P. Mehendale }
\email{dhananjay.p.mehendale@gmail.com}
\affiliation{Sir Parashurambhau College,Pune,India-411007.}

\date{\today}
 
\begin{abstract} 

In this paper we propose an algorithm to generate all possible graceful graphs (including trees) containing $n$ vertices 
as lattice paths in certain triangular lattice defined below. This lattice that corresponds to graphs containing 
$n$ vertices is called $n$-lattice and is made up of certain rows of vertex pairs $(i,j).$ Each row 
of this $n$-lattice is made up of those vertex-pairs, say $(i,j),$ for which $|i-j|$ is same for every vertex-pair 
belonging to that row, and where $i,j\in \{1,2,\cdots, n\}.$ The first row of this $n$-lattice contains $(n-1)$ 
vertex-pairs, $(i,i+1), i = 1,2,\cdots,(n-1).$ The second row of this $n$-lattice contains $(n-2)$ 
vertex-pairs, $(i,i+2), i = 1,2,\cdots,(n-2).$ In this way, one goes down to the last row of this lattice which 
contains only one vertex-pair, $(1,n).$ 
A lattice path is the one made up of $(n-1)$ vertex-pairs such that every row of the triangular lattice contributes 
exactly one vertex-pair to this lattice path. We obtain all possible lattice paths without omission or repetition by 
generating them in a systematic way in a well-defined lexicographic order. The collection of all such lattice paths 
forms all possible graceful graphs. We will note various observations related to these lattice paths. For example,  
the lattice paths appear in symmetric pairs, i.e. for each lattice path there exists a corresponding unique lattice path 
which is the mirror images of this lattice path taken in the line of symmetry passing vertically and centrally through 
the lattice and each lattice path and its corresponding mirror image represent isomorphic graceful graphs. The main result 
of this paper is the affirmative settlement of the well-known graceful tree conjecture. 

\end{abstract}
 

\maketitle

\section{Introduction:}

A tree on n vertices is said to be graceful or said to have a graceful labeling if when its vertices are labeled with 
numbers $\{1, 2, ...,n\}$ and edges are labeled by the difference, $|i-j|,$ of their respective end vertex 
labels, $i,j \in \{1, 2, ...,n\}$ then all the edge labels taken together constitute the set $\{1,2, ..., n-1\}.$ 

In the year 1964 Ringel \cite{Ringel} proposed the following conjecture:

{\bf Conjecture 1.1(Ringel):} If $T$ is a fixed tree with $m$ edges, then the complete graph, $K_{(2m+1)},$ on $(2m+1)$ 
vertices, can be decomposed into $(2m+1)$ copies of $T.$

Attempts to prove Ringel’s conjecture have focused on the following stronger conjecture about trees, the so called  
Graceful Tree Conjecture \cite{West}: 

{\bf Conjecture 1.2 (Graceful Tree Conjecture):} A tree of each isomorphic type is graceful, i.e. has a graceful labeling. 

As an implication of the validity of the Graceful Tree Conjecture for a particular tree $T$ under consideration 
Rosa \cite{Rosa} obtained the following result: 

{\bf Theorem 1.1(Rosa):} If a tree $T$ with $m$ edges has a graceful labeling then $K_{(2m+1)}$ has decomposition 
into $(2m+1)$ copies of $T.$

Thus, the validity of Ringel’s conjecture will automatically follow from the validity of the graceful tree conjecture. 

In this paper we settle Graceful Tree Conjecture in the affirmative by showing the existence of graceful labeling for 
a tree of each isomorphism type. We give a systematic procedure to generate in a well defined lexicographic order all the 
graceful graphs containing $n$ vertices that exist as subgraphs or can be constructed in a systematic way without omission 
or repetition as subgraphs of a complete graph on $n$ vertices which have been labeled by numbers $\{1, 2, ...,n\}.$  

Is graceful labeling unique? No. There can exist several distinct graceful labellings for a given unlabeled graph. 
The following simple result provides two distinct graceful labellings for the same graph, i.e. graph of same isomorphism type.

{\bf Theorem 1.2:} Every graceful $(n,n-1)$ tree remains graceful under the
mapping of vertex labels $f:\{1,2,\cdots,n\}\to \{1,2,\cdots,n\}$ such that $j\to(n-j+1),$ i.e. $f(j)=(n-j+1).$

{\bf Proof:} It is easy to check that under this mapping $1\to(n-1),2\to(n-2),3\to(n-3),\cdots,n\to1.$ Thus, this mapping is 
one-to-one and onto. It is easy to check that the edge represented by the vertex-pair $(i,k)$ changes to the edge 
represented by the vertex-pair $(n-i+1,n-k+1)$ under this map. Also, note that the edge label for the 
edge represented by the vertex-pair $(i,k)$ will be $|i-k|$ and also under the above mentioned mapping the edge label 
of the new changed edge represented by the vertex-pair $(n-i+1,n-k+1)$ will also 
be $|(n-i+1)-(n-k+1)|=|i-k|.$ Thus, after replacing the vertex labels, $j,$ of a graceful $(n,n-1)$-tree 
by the vertex labels $(n-j+1)$ the graceful tree under consideration will remain graceful as before under the 
said mapping.

Further we see below that every graceful graph (which may not be a tree) also remains graceful under the above given 
same mapping.

{\bf Theorem 1.3:} Every graceful graph containing $n$ vertices, $\{1,2,\cdots,n\}$ remains graceful under the
mapping of vertex labels: $j\to(n-j+1)$

{\bf Proof:} This mapping is clearly one-to-one and onto as seen above. Let $i,k$ be the vertex labels of two adjacent 
vertices of the given graceful graph. Then the edge label for this edge will be $|i-k|.$ Now under the above 
mentioned mapping the edge label still remains the same because when we change the vertex labels $i,k$ to new vertex labels 
$(n-i+1,n-k+1)$ under the specified mapping the edge label does not change because the edge label   
$|i-k|$ changes to $|(n-i+1)-(n-k+1)|$ which is same as $|i-k|.$ Thus, after replacing the vertex labels, $j,$ of given 
graceful graph by the vertex labels $(n-j+1)$ the graceful graph under consideration still remains graceful. 

In the next theorem we will determine the total number of graceful graphs on $n$ vertices that are possible.

{\bf Theorem 1.4:} There exist in all $(n-1)!$ graceful graphs on $n$ vertices. 

{\bf Proof:} We form all possible graceful graphs on $n$ vertices as spanning subgraphs of complete graph on $n$ 
vertices, $K_n,$ whose vertices we label by integers $\{1, 2, ...,n\}.$ We form the following spanning subgraphs, one by one, 
of this $K_n$ whose vertices are labeled by integers $\{1, 2, ...,n\}$ such that they are made up of all possible distinct 
sets of edges. We represent these edges as vertex-pairs $(i,j)$ where $i,j$ are labels of the end vertices of that edge. 
Each such set is made up of exactly $(n-1)$ edges like: 

$\{(i_{(n-1)},j_{(n-1)}),(i_{(n-2)},j_{(n-2)}),\cdots,(i_1,j_1)\}$ 

such that the edges chosen to form this set satisfy the following conditions:

$(i_{(n-1)},j_{(n-1)}) \in \{(1,n)\}$

$(i_{(n-2)},j_{(n-2)}) \in \{(1,n-1),(2, n)\}$

$(i_{(n-3)},j_{(n-3)}) \in \{(1,n-2),(2, n-1),(3, n)\}$

$\vdots$

$(i_1,j_1) \in \{(1,2),(2,3),(3,4),\cdots,(n-1,n)\}$ 

and view it as a subgraph of $K_n.$ 

Thus, the first edge represented by the vertex-pair $(i_{(n-1)},j_{(n-1)})$ has only one choice, namely, $(1,n)$ and 
The second edge represented by the vertex-pair $(i_{(n-2)},j_{(n-2)})$ has two choices, $\{(1,n-1),(2, n)\}.$ 
By proceeding in this way the last edge represented by the vertex-pair $(i_1,j_1)$ has $(n-1)$ choices   
$\{(1,2),(2,3),(3,4),\cdots,(n-1,n)\}.$ Note that with these choices 
all the labeled graphs thus formed are automatically graceful by choice and all these choices 
can be made independently. Therefore the total number of graceful graphs one can form containing $n$ vertices 
will be equal to $1\ts2\ts3\ts\cdots\ts(n-1) = (n-1)!$

\section{Generating $(n+1,n)$-trees from $(n,n-1)$-trees:} 

A graph, $G,$ with $p$ vertices and $q$ edges will be called a $(p,q)$-graph, $G.$ We denote the vertex set 
(set of vertices) for this graph by the set $V(G)=\{v_1,v_2,\cdots,v_p\}$ and the edge set (set of edges) for this 
graph, $E(G)=\{e_1,e_2,\cdots,e_q\},$ where $e_k=(v_i,v_j)$ where $v_i,v_j \in V(G)$ are the end vertices of the 
edge $e_k\in E(G).$
A graph is called connected if there exists at least one path joining every vertex to every other vertex. For a tree 
graph there exists exactly one path joining any two vertices. A tree with $n$ vertices is a connected $(n,n-1)$-graph. 
A graph is called acyclic if there does not exists a closed loop made up of its edges. 
A tree with $n$ vertices is an acyclic $(n,n-1)$-graph. 

{\bf Definition 2.1:} Let $T$ be an $(n,n-1)$-tree of some isomorphism type. Let $T^e$ be a tree obtained by taking 
a (new) vertex $u$ outside of the set of vertices $V(T)$ of $T$, i.e. $u\not\in V(T),$ and joining $u$ to any one 
vertex $v\in V(T)$ by an (new) edge represented as vertex-pair $(u,v)$ not in the set of edges, $E(T),$ of $T.$ 
In this case the tree $T^e=T+(u,v)$ is called an {\bf extension} of the tree $T$ (obtained by extending $T$ at $v\in V(T)$).    

{\bf Definition 2.2:} The subset $V_j$ of vertices $\{u_1,u_2 ,\cdots,u_r\}$ of the vertex set $V(T)$ of tree $T$ 
is called a set of equivalent vertices or simply an {\bf equivalent set} if all the trees $T' = T+(v,u_s),1\le s\le r,$ 
obtained by extending the tree $T$ by adding an edge $(v,u_s),$ obtained by joining vertex $u_s\in V_j$ to a new 
vertex $v\not\in V(T),$ are all isomorphic to each other. 

{\bf Definition 2.3:} The collection of subsets $\{V_1,V_2,\cdots,V_m\}$ of $V(T),$ the
vertex set of tree $T$ is called a {\bf partitioning} of $V(T)$ into equivalent
sets if all the subsets $V_i,i=1,2,\cdots,m$ are all possible different equivalent sets, such that $V_i\cap V_j=\varnothing$ 
when $i\neq j$ where $\varnothing$ is an empty set, and $V(T)=\cup_{i=1}^mV_i.$ 

{\bf Definition 2.4:} The vertices $w_i$ and $w_j$ are called {\bf inequivalent vertices} if they belong to different 
equivalent sets, i.e if $w_i\in V_i$ while $w_j\in V_j$ and $V_i,V_j$ are different equivalent sets having null intersection.  

{\bf Definition 2.5:} The {\bf $(n,n-1)$-stock} is the collection of all possible unlabeled (non-isomorphic)
$(n,n-1)$-trees.

{\bf The Systematic Extension Procedure of Trees in $(n,n-1)$-stock:} We now discuss in simple steps the systematic 
extension procedure to be applied to the trees in the $(n,n-1)$-stock to create $(n+1,n)$-stock:

(1) Take a tree $T$ in $(n,n-1)$-stock.

(2) Perform partitioning of the vertex set $V(T)$ of tree $T$ into all possible disjoint subsets such that each of this 
subset is an equivalent set, i.e. 
$V(T)=\cup_{i=1}^mV_i$ such that $V_i,i=1,2,\cdots,m$ are subsets which are equivalent sets, 
such that $V_i\cap V_j=\varnothing,$ where $\varnothing$ is an empty set, when $i\neq j.$ 

(3) Choose some (representative) vertex $u_j\in V_i,i=1,2,\cdots,m$ from each equivalent set in the partitioning and 
extend $T$ by joining $u_j$ to a new vertex $v_i\not\in V(T)$ i.e. obtain $T^j_i = T + (u_j,v_i).$ The new trees $T^j_i$ 
are essentially all the non-isomorphic trees that are possible to obtain from tree $T$ of some isomorphism type belonging 
to $(n,n-1)$-stock.

(4) Repeat the above steps (1)-(3) for every tree in $(n.n-1)$-stock.

Steps (1)-(4) depicts the systematic extension procedure applied on the trees in the $(n,n-1)$-stock to produce the trees 
in the $(n+1,n)$-stock.

{\bf Theorem 2.1:} The application of the systematic extension procedure on each and every tree of the trees in the 
$(n,n-1)$-stock produces $(n+1,n)$-stock, i.e. the $(n+1,n)$-stock emerges from 
the $(n,n-1)$-stock by applying the systematic extension procedure described above on all the trees in the $(n,n-1)$-stock. 

{\bf Proof:} Every tree $T^*$ belonging to $(n+1,n)$-stock can be considered as arrived at by extension at some vertex $u_j$ 
belonging to some equivalent set of vertices which is a subset of $V(T)$ of a tree $T$ belonging to $(n,n-1)$-stock. 

\section{Lattice, Lattice Paths, and Graceful Labeling:} 

We assign as the vertex labels to any graph containing $n$ vertices the integers $\{1,2,3,\cdots,n\},$ 
i.e. $V(G) = \{1,2,3,\cdots,n\}.$ We represent edges of the same graph in terms of the vertex-pairs $(i,j)$ of adjacent 
vertices where $i,j$ are labels of the end vertices of the corresponding edge represented as a vertex-pair $(i,j).$ 
Thus, each vertex-pair $(i,j)$ corresponds to an edge whose end vertices are labeled $i$ and $j$ respectively. 

{\bf Definition 3.1:} The {\bf $n$-lattice} is a triangular pattern, having shape of an inverted triangle, containing rows 
of vertex-pairs, $(i,j),$ representing corresponding edges with end vertices $i,j$ as described below:

In the first (or the top row) of this so called $n$-lattice there are vertex-pairs 

$(1,2),(2,3),(3,4),\cdots,(n-1,n).$

In the second (or the row just below the top row) there are vertex-pairs 

$(1,3),(2,4),(3,5),\cdots,(n-2,n).$

In the third (or the row just below the second row) there are vertex-pairs 

$(1,4),(2,5),(3,6),\cdots,(n-3,n).$

Continuing on these lines, in the second last or $(n-2)$-th row there are vertex-pairs

$(1,n-1),(2,n)$

and in the last or $(n-1)$-th row there is only one vertex-pair 

$(1,n).$ 

In other words, the $n$-lattice is a triangular shaped lattice having shape of an inverted triangle, 
such that certain distinct vertex-pairs (each representing a unique edge) are written in rows as follows: 
The first or top row contains vertex-pairs $(i,i+1),$ where $i$ goes from 1 to $n-1.$ 
The second row written below the first row contains vertex-pairs $(i,i+2),$
where $i$ goes from 1 to $n-2.$ The $k$-th row, reached
by successively creating rows downwardss, contains vertex-pairs $(i,i+k),$
where $i$ goes from 1 to $n-k.$ The last row contains a single vertex-pair $(1,n).$  
Thus, the so called $n$-lattice is the one as shown below:
   \begin{align*}
    (1,2)(2,3)(3,4)(4,5)(5,6)\cdots(n-1,n)     
   \end{align*}
   \begin{align*}
   (1,3)(2,4)(3,5)(4,6)\cdots(n-2,n)
   \end{align*}
   \begin{align*}
   (1,4)(2,5)(3,6)\cdots (n-3,n)
   \end{align*}
   \begin{align*}
   \vdots 
   \end{align*}
   \begin{align*}
   (1,n-2)(2,n-1)(3,n)
   \end{align*}
   \begin{align*}
   (1,n-1)(2,n) 
   \end{align*}
   \begin{align*}
   (1,n)
   \end{align*}
   
{\bf Definition 3.2:}If we choose first the vertex-pair from the last (bottom) row, i.e. $(1,n),$ then from the a vertex-pair 
in the second-last row, i.e. a vertex-pair $\in \{(1,n-1),(2,n)\},$ then from the third-last row and so on, $\cdots,$ 
and finally from the first row containing vertex-pairs $(1,2),(2,3),(3,4),(4,5),(5,6),\cdots,(n-1,n)$ then such a collection 
of vertex-pairs is called a lattice path belonging to the $n-$lattice. Since 
firstly we have chosen the vertex-pair $(1,n)$ belonging to the last row of the $n-$lattice then we have chosen a vertex-pair 
either $(1,n-1)$ or $(2,n)$ in the second-last row and so on, therefore, this lattice path can be further distinguished 
by calling it as a {\bf upwards lattice path} in the $n$-lattice. We represent this lattice path by writing the vertex-pairs, 
one after the other, by first writing the chosen vertex-pair from the last row then from the 
second last row and so on till we finally write the vertex-pair chosen from the first row. 

These (upwards) lattice paths can be formed in a lexicographic order so that there will be neither any repetition nor any 
omission if we form them in an orderly way as follows: 

$(1,n)(1,n-1)(1,n-2)\cdots (1,3)(1,2)$ $\to$ smallest in the lexicographic order

$(1,n)(1,n-1)(1,n-2)\cdots (1,3)(2,3)$

$(1,n)(1,n-1)(1,n-2)\cdots (1,3)(3,4)$

$\cdots$

$\cdots$

$\cdots$

$(1,n)(1,n-1)(1,n-2)\cdots (1,3)(n-1,n)$

$\cdots$

$\cdots$

$\cdots$

$(1,n)(2,n)(3,n)\cdots (n-1,n)$ $\to$ largest in the lexicographic order
   
{\bf Definition 3.3:} If we choose first the vertex-pair from the first (top) row, then on the second row, then on the 
third row, $\cdots,$ and finally on the $(n-1)$-th (last) 
row then such a collection of vertex-pairs is called a lattice path belonging to the $n-$ lattice. Since firstly we 
have chosen the vertex-pair belonging to first (top) row, then from the second row, then from the third row, $\cdots,$ 
and finally from the $(n-1)$-th (last) row such a lattice path can be further distinguished 
by calling it as a {\bf downwards lattice path} in 
the $n$-lattice. We represent this lattice path by writing the vertex-pairs in a sequence 
starting with writing firstly the chosen vertex-pair from the first row then writing the vertex-pair chosen from the 
second row and so on till we finally write the vertex-pair chosen from the last ($(n-1)$-th) row. 

These (downwards) lattice paths can be formed in a lexicographic order so that there will be neither any repetition nor any 
omission if we form them in an orderly fashion as follows: 

$(1,2)(1,3)(1,4)\cdots (1,n-2)(1,n-1)(1,n)$ $\to$ smallest in the lexicographic order

$(1,2)(1,3)(1,n-2)\cdots (1,n-2)(2,n)(1,n)$

$\cdots$

$\cdots$

$\cdots$

$(n-1,n)(1,3)(1,4)\cdots (1,n-2)(1,n-1)(1,n)$

$\cdots$

$\cdots$

$\cdots$

$(n-1,n)(n-2,n)(n-3,n)\cdots (2,n)(1,n)$ $\to$ largest in the lexicographic order

{\bf Example:}

The 4-lattice is as as given below:

\begin{align*}
    (1,2) (2,3)(3,4)     
   \end{align*}
   \begin{align*}
   (1,3)(2,4)
   \end{align*}
   \begin{align*}
   (1,4)
   \end{align*}
   
The upwards lattice paths in the above 4-lattice are 

$(1,4)(1,3)(1,2)$ $\to$ smallest in the lexicographic order

$(1,4)(1,3)(2,3)$

$(1,4)(1,3)(3,4)$

$(1,4)(2,4)(1,2)$

$(1,4)(2,4)(2,3)$

$(1,4)(2,4)(3,4)$ $\to$ largest in the lexicographic order

The downwards lattice paths in the above 4-lattice are 

$(1,2)(1,3)(1,4)$ $\to$ smallest in the lexicographic order

$(1,2)(2,4)(1,4)$

$(2,3)(1,3)(1,4)$

$(2,3)(2,4)(1,4)$

$(3,4)(1,3)(1,4)$

$(3,4)(2,4)(1,4)$ $\to$ largest in the lexicographic order.

Thus, a lattice path (of upwards type) is a path obtained by selecting some one
vertex-pair from each row of the $n$-lattice and writing down these vertex-pairs 
in sequence starting with the vertex-pair $(1,n)$ on the lowest row (i.e. the row at the 
bottom of the lattice containing only one vertex-pair, $(1,n)$) and moving
up in succession incorporating the chosen vertex-pair on each row till we finally select a 
vertex-pair on the top (or the first) row of the $n$-lattice and the lattice 
path is thus formed. 

On the other hand, a lattice path (of downwards type) is a path obtained by selecting some one
vertex-pair from each row of the $n$-lattice and writing down these selected vertex-pairs 
in sequence starting with a vertex-pair on the first (top) row and moving
down in succession incorporating the chosen vertex-pair on each row till the lattice 
path is finally formed with the selection of vertex-pair $(1,n)$ on the last row at the bottom of the lattice.

Existence of a lattice path in the $n$-lattice that corresponds to a tree of some isomorphism type then it 
implies the existence of a graceful labeling for the tree of that isomorphism type. 
Note that in a lattice path every row of the lattice contributes exactly one vertex-pair  
therefore every lattice path essentially corresponds to a graceful graph of some isomorphism type. In other words, 
since by definition a lattice path contains exactly one edge of type $(i,j)$ such 
that $|i-j|=k$ for every $k\in\{1,2,3,\cdots,n-1\},$ therefore, every lattice path certainly represents a graceful 
graph of some isomorphism type.

Note that a lattice path in the $n$-lattice corresponds to a graceful tree when it is an $(n,n-1)$-graph which is a 
connected graph or an $(n,n-1)$-graph which is an acyclic graph. We are mainly interested in this paper to see whether 
every tree, i.e. a tree of each isomorphism type is graceful, i.e. has a graceful labeling.  

In order to get the pictorial representation of a lattice path, in the triangular lattice shown above we join by a 
line segment the vertex-pair $(1,n)$ chosen from last row at the bottom of the 
triangular lattice to the vertex-pair $(1,n-1)$ or $(2,n)$ that has been chosen in the second last row and then join 
this vertex-pair chosen from the second last row to the vertex-pair that has been chosen in the third last row 
and so on till we finally join the vertex-pair that has been chosen in the second to the vertex-pair that has been chosen 
in the first row, and thus complete the formation of the the pictorial representation of the lattice path under consideration. 

Alternatively, if we will join the vertex-pair chosen from first row to the vertex-pair that has been chosen in 
the second row and then the vertex-pair that has been chosen from second row to the vertex-pair that has been chosen 
in the third row and so on till we finally join the vertex-pair that has been chosen in the second last 
row ($(n-2)$-th row) to the vertex-pair that has been chosen in the last 
row ($(n-1)$-th row) then this also leads to the pictorial representation of the lattice path under consideration.

Note that since each lattice path contains exactly one 
vertex-pair from each row of lattice, $(i,j)$ say, having the associated edge label $|i-j|=k$ when the 
chosen vertex-pair belongs to $k$-th row this clearly implies that each lattice path corresponds to a graceful graph. 

Further note that if we will replace each vertex label say $j$ appearing in the vertex-pair belonging to a 
lattice path by the vertex label $(n-j+1)$ we will get a new lattice path corresponding to a graceful graph which is 
pictorially the mirror image taken in the vertical axis passing vertically 
upwards through vertex-pair $(1,n)$ of the earlier obtained lattice path (Refer to Theorems 1.2 and 1.3 above). 
More formally, 

{\bf Definition 3.4:} An imaginary vertical line starting from vertex-pair $(1,n)$ and going upwards passing through 
the vertex-pairs $(2,n-1),(3,n-2),\cdots,$ rising up to first row is called the {\bf line of symmetry} for the $n$-lattice. 

It is easy to see that the line of symmetry is dividing the $n$-lattice into two equal parts.

The pictorial representation of a given lattice path and the pictorial representation of the lattice path obtained from 
this given lattice path by replacing each vertex label, say $j,$ appearing in the vertex-pairs of this lattice path 
by the vertex label $(n-j+1)$ are {\bf mirror image} of each other taken in the line of symmetry. It is easy to check 
that the graceful graphs corresponding to the pictorial representations lattice paths that are mirror images of each other 
are isomorphic, i.e. they form the two gracefully labeled copies of the same (unlabeled) graph up to isomorphism. For example 
the pictorial representations of the following two lattice paths  
$(1,n)(1,n-1)(1,n-2)\cdots(1,2)$ and $(1,n)(2,n)(3,n)\cdots(n-1,n)$ are mirror images of each other in the line of symmetry. 
It is easy to check that these lattice paths are lying at left 
and right boundary of the $n$-lattice and correspond to gracefully labeled $(n,n-1)$-trees which are star trees. 
It is easy to check that the following two zigzag pictorial representations of the lattice paths, 
namely, those obtained by joining vertex-pairs $(1,n)(1,n-1)(2,n-1)(2,n-2)(3,n-2)\cdots$
and $(1,n)(2,n)(2,n-1)(3,n-1)(3,n-2)(4,n-2)\cdots$ are mirror images of each other around the line of symmetry. 
It is easy to check that these zigzag lattice paths correspond to gracefully labeled $(n,n-1)$-trees which are path-trees. 

Note that since the upwards lattice paths and the downwards lattice paths are essentially same in content so hereafter 
we will be considering upwards lattice paths only. So, hereafter, {\bf a lattice path meant an upwards lattice path,} i.e. 
hereafter all lattice paths are upwards lattice paths.

\section{Some Related Observations:}

{\bf Observation 4.1:} From a given graceful tree one can obtain another graceful tree with different labeling isomorphic 
to the original tree by taking mirror image of the pictorial representation of the associated lattice path taken in the 
line of symmetry.  

{\bf Theorem 4.1:} If there exists graceful labeling to an $(n,n-1)$-tree, $T$, of some isomorphism type then there exist 
at least two graceful copies with different labeling for that same tree.  

{\bf Proof:} Since $(n,n-1)$-tree, $T$, has graceful labeling therefore there exists a lattice path say, 
$(1,n)(2,n)\cdots (p,p+1)$ where $p \in \{1,2,\cdots,(n-1)\}.$ Construct the pictorial representation of this lattice 
path in $n$-lattice by joining the vertex-pairs by line segments by starting the vertex-pair $(1,n)$ in the bottom most 
row and joining it to the vertex-pair $(2,n)$ in the row just above it by a line segment and so on 
in the same order. Construct mirror image of this lattice path in the line of symmetry of the $n$-lattice. This mirror image 
will be the lattice path $(1,n)(1,n-1)\cdots (n-p,n-p+1)$ and construct as above its pictorial representation. These two 
pictorial representations are mirror images of each other taken in the line of symmetry. They represents 
graceful trees isomorphic to each other with different labeling.  

{\bf Example:}

Consider a 6-lattice as given below:

\begin{align*}
    (1,2)(2,3)(3,4)(4,5)(5,6)     
   \end{align*}
   \begin{align*}
   (1,3)(2,4)(3,5)(4,6)
   \end{align*}
   \begin{align*}
   (1,4)(2,5)(3,6)
   \end{align*}
   \begin{align*}
   (1,5)(2,6)
   \end{align*}
   \begin{align*}
   (1,6)
   \end{align*}
   
   Consider a lattice path: 
   
   $(1,6)(1,5)(2,5)(1,3)(3,4)$ 
   
   Its mirror image in the vertical line called the line of symmetry (passing through vertex-pairs $(1,6)(2,5)(3,4)$) will be 
   
   $(1,6)(2,6)(2,5)(4,6)(3,4)$.
   
   Note that actually this mirror image is the same one obtained by map 
   
   $\Theta:V(T)\to V(T)$ defined by $j\to (n+1)-j,$ with $n=6$ in the 
   above case. 
   
{\bf Observation 4.2:} Apart from mirror images there exist multiple graceful copies with different labellings for 
the tree of same isomorphism type.   

For example, trees represented by lattice paths 

$(1,n)(1,n-1)(1,n-3)\cdots (1,3)(2,3)$ and 

$(1,n)(1,n-1)(2,n-1)(3,n-1)\cdots (n-3,n-1)(n-2,n-1)$ 

are isomorphic though they are not mirror images of each other taken in the line of symmetry. 
Though not mirror images taken in the line of symmetry, these trees have same shape as a lattice 
paths representing a tree which is juxtaposition (join) of an $(n-2)$-star tree and a 1-path-tree.

{\bf Definition 4.1:} A lattice path is called a partial (or incomplete) lattice path if it contains vertex-pairs 
belonging some $k$ rows, $k<n,$ under consideration, i.e. it is a part of (an entire or complete) lattice path 
which contains exactly one vertex-pair from each row of the $n$-lattice. 
   
{\bf Theorem 4.2:} The count of graceful $(n,n-1)$ trees in an $n$-lattice is $\ge 2^{(n-2)}.$ 

{\bf Proof:} We construct lattice paths as follows: 

(i) We start with the vertex-pair $(1,n).$ 

(ii) We append to it a vertex-pair in the row just above it leading to two new extended partial lattice paths

$(1,n)(1,n-1)$ and $(1,n)(2,n).$

(iii) We continue appending each time a vertex-pair in the row just above it which adds only one new vertex through that 
appended vertex-pair by the following 

{\bf Rule:} Append each time a vertex-pair in the row just above such that one label of this appended vertex-pair is 
same as one of the vertex present the vertex-pair appended in the earlier stage and the other label in this appended 
vertex-pair differs from the vertex in the earlier appended vertex-pair by plus or minus one. 

This rule keeps the resulted partial lattice path a tree. In short, each of the partial lattice paths grows in 
two ways leading each time to two new partial (incomplete) lattice paths which are trees till we reach the top row where 
we get a (full) lattice path (which is no more partial) and the construction of graceful tree completes.    

(iv) Thus, we continue appending suitable vertex-pairs satisfying the above simple rule leading to two new 
partial (incomplete) lattice-paths from each of the partial (incomplete) lattice-paths obtained in the earlier stage.
Thus, each of the above lattice paths grow in two ways in accordance with the above rule leading to two new partial 
(incomplete) lattice paths from each of them as given below. Thus we get two lattice paths from each lattice path arrived 
at the earlier stage by extending lattice paths of the earlier stage as per the above rule as follows: 

We get partial lattice paths $(1,n)(1,n-1)(1,n-2)$ and $(1,n)(1,n-1)(2,n-1)$ from the partial lattice path $(1,n)(1,n-1)$ and 

We get partial lattice paths $(1,n)(2,n)(2,n-1)$ and $(1,n)(2,n)(3,n)$ from the partial lattice path $(1,n)(2,n).$ 

At each stage of moving to upper row and appending new vertex-pairs it is clear to see that each earlier partial lattice 
path gives rise to two new partial lattice paths. Thus, at each step of moving to upper row and appending new vertex-pairs 
the count of partial lattice paths becomes double. If we continue growing these partial lattice paths by appending new 
vertex-pairs in the upper row following the rule then finally we reach the stage of appending the vertex-pairs from the 
first row and the lattice paths thus arrived at are no more partial and they become complete and it is straightforward 
to check that they all now correspond to graceful trees. 

Let $n = 2.$ In 2-lattice there exists only one, $2^{(n-2)} = 2^0 = 1,$ graceful tree represented by lattice 
path $(1,2),$ a one edge graceful tree. For this case, there are $n-2 = 0$ stages of appending giving rise 
to $2^{(n-2)} = 2^0 = 1,$ graceful tree.

Let $n = 3.$ The 3-lattice is as given below:

\begin{align*}
    (1,2) (2,3)     
   \end{align*}
   \begin{align*}
   (1,3)
   \end{align*}
   
In this lattice there exist two, $2^{(n-2)} = 2^1 = 2,$ graceful trees represented by lattice paths 

$(1,3)(1,2)$ and $(1,3)(2,3).$ For this case, there are $n-2 = 1$ stages of appending giving rise 
to $2^{(n-2)} = 2^1 = 2,$ graceful trees.

Let $n = 4.$ The 4-lattice is as given below:

\begin{align*}
    (1,2) (2,3)(3,4)     
   \end{align*}
   \begin{align*}
   (1,3)(2,4)
   \end{align*}
   \begin{align*}
   (1,4)
   \end{align*}
   
The four, $2^{(n-2)} = 2^2 = 4,$ graceful trees in this lattice are represented by the following lattice paths:

$(1,4)(1,3)(1,2)$ 

$(1,4)(1,3)(2,3)$

$(1,4)(2,4)(2,3)$

$(1,4)(2,4)(3,4).$ 

For this case, there are $n-2 = 2$ stages of appending giving rise 
to $2^{(n-2)} = 2^2 = 4,$ graceful trees. 

Continuing on these lines one can check that for the case of $n$-lattice there will be $n-2$ stages of appending giving rise 
to $2^{(n-2)}$ graceful trees obtained by observing rule defined in (ii). It is easy to check further that if we give some 
relaxation in the above rule and if we allow to append new vertex-pairs in the next upper row of $n$-lattice such that the 
newly appended vertex-pair should add one (or more) new vertices at each appending stage then we can easily get more new 
graceful trees not belonging to the set of newly formed graceful trees following above rule.  
Hence the result.
   
{\bf Observation 4.3:} Let $T$ be an $(n,n-1)$-tree  of some isomorphism type. Let the vertex set, $V(T),$ of tree $T$ be 
partitioned into equivalent sets. It is clear to see that for every vertex of the tree $T$ there exists some unique 
equivalent set to which it belongs. Suppose there exist in all $m$ equivalent 
sets such that $V(T)=\cup_{i=1}^mV_i,$ and $V_i\cap V_j=\varnothing$ when $i\neq j$ 
where $\varnothing$ is an empty set. Now, by extending at some one vertex in each such equivalent set will lead to  
$m$ new (non-isomorphic) $(n+1,n)$-trees obtained from this $(n,n-1)$-tree $T.$ 

{\bf Theorem 4.3:} The cardinality of new (mutually non-isomorphic) $(n+1,n)$-trees that can emerge from a 
given $(n,n-1)$-tree, $T,$ by extension (see definition 2.1) is equal to the cardinality the equivalent sets 
associated with tree $T.$ 

{\bf Proof:} Each equivalent set of vertices corresponding to given given $(n,n-1)$-tree, $T,$ gives rise to exactly 
one and only one new (mutually non-isomorphic) $(n+1,n)$-tree through extension by attaching a new edge joining any 
one vertex belonging to the equivalent set under consideration to a new vertex not belonging to the vertex set, $V(T),$ 
of the given tree, $T.$

{\bf Observation 4.4:} The $n$-lattice contains in terms of lattice paths sufficiently many copies of graceful 
$(n,n-1)$-trees of each isomorphism type so that label ``1'' will be automatically present at a vertex belonging to some   
equivalent set in some gracefully labeled copy of that $(n,n-1)$-tree, $T,$ of some isomorphism type such that when this 
gracefully labeled tree will be extended by attaching new edge represented by vertex-pair $(1,n+1)$ we get a graceful 
copy of the desired $(n+1,n)$-tree. Thus, we can obtain graceful tree corresponding to any unlabeled $(n+1,n)$-tree 
(i.e. an $(n+1,n)$-tree of each isomorphism type) by attaching new edge represented by vertex-pair $(1,n+1)$ at some 
vertex belonging to some equivalent set of some graceful copy of an $(n,n-1)$-tree of some isomorphism type. 

{\bf Theorem 4.4:} For any given unlabeled $(n+1,n)$-tree there exists a suitable gracefully labeled $(n,n-1)$-tree, $T$ say 
of some isomorphism type having a vertex with label ``1'' belonging to some of its equivalent sets such that when $T$ is 
extended by appending a new vertex-pair $(1,n+1)$ we get the gracefully labeled tree $T^e = T + (1,n+1)$ which is isomorphic 
to the given unlabeled $(n+1,n)$-tree under consideration. Further, this statement is true for every unlabeled $(n+1,n)$-tree. 

{\bf Observation 4.5:} For every unlabeled $(n+1,n)$-tree $T^*$ there exists at least one graceful 
$(n,n-1)$-tree $T$ of some isomorphism type having a vertex with label ``1'' in some of its equivalent set such that 
we get gracefully labeled version of tree $T^*$ by extending the graceful tree $T$ at vertex with label '1' by appending edge 
$(1,n+1)$ i.e. $T^*$ is isomorphic to $T+(1,n+1).$ 

{\bf Theorem 4.5:} For every unlabeled $(n+1,n)$-tree $T^*$ there exists at least 
one $(n,n-1)$-tree $T,$ of some isomorphism type such that a graceful version among the graceful versions of this  
$(n,n-1)$-tree, $T$ say, has a vertex (in some equivalent set) with label ``1''  such that we get gracefully labeled 
version of tree $T^*$ by extending this graceful version, $T,$ at this vertex with label '1' by appending edge 
$(1,n+1)$ i.e. $T^*$ is isomorphic to gracefully labeled $T$ appended by $(1,n+1)$ i.e. $T^*$ is isomorphic to $T+(1,n+1),$ 
and therefore there exists a graceful labeling for tree $T^*$ i.e. $T^*$ is graceful. 

{\bf Example:} There are three non-isomorphic trees containing {bf five vertices} (i) 4 star, (ii) (3 star-1 path-tree) 
(iii) 4 path-tree. 

Let us consider the following six lattice paths in the $5$-lattice representing graceful versions for above trees 
containing five vertices: 

(1) $(1,5)(1,4)(1,3)(1,2)$ 

(2) $(1,5)(1,4)(1,3)(2,3)$

(3) $(1,5)(1,4)(2,4)(2,3)$

(4) $(1,5)(1,4)(2,4)(3,4)$

(5) $(1,5)(1,4)(3,5)(1,2)$

(6) $(1,5)(1,4)(3,5)(2,3)$

and further there are six more lattice paths which are mirror images of above lattice paths taken in the line of symmetry:  

(7) $(1,5)(2,5)(1,3)(3,4)$ 

(8) $(1,5)(2,5)(1,3)(4,5)$

(9) $(1,5)(2,5)(2,4)(2,3)$

(10) $(1,5)(2,5)(2,4)(3,4)$

(11) $(1,5)(2,5)(3,5)(3,4)$

(12) $(1,5)(2,5)(3,5)(4,5)$ 

We now see below how vertex with label ``1'' is present at a vertex in sufficiently many equivalent sets of 
sufficiently many gracefully labeled non-isomorphic trees such 
that by extending these labeled trees by appending them by vertex-pair $(1,6)$ we get graceful trees for trees 
of each isomorphism type containing six vertices.

{\bf (a)} The graceful trees given by lattice paths (1) and (12) are 4 star trees. We give below the equivalent sets for 
these trees in (1) and (12):

For tree (1) the partitioning of vertices into equivalent sets is $\{1\},\{2,3,4,5\}$. Here label ``1'' is present in 
the first equivalent set. 

For tree (12) the partitioning of vertices into equivalent sets is $\{5\},\{1,2,3,4\}$. Here label ``1'' is present in 
the second equivalent set. 

If we will now extend the graceful trees (1) and (12) by appending vertex-pair $(1,6)$ we will get graceful trees, 
respectively a (5 star) and a (4 star - 1 path-tree), and these are the only possible non-isomorphic trees one can 
get by extending a (4 star) tree. 

{\bf (b)} The graceful trees given by lattice paths (2),(4),(5),(8),(9) and (11) are all isomorphic, are (3 star-1 path-tree) 
type trees. 

If we will draw corresponding {\bf tree diagrams} from the above lattice paths, (2),(4),(5),(8),(9) and (11), we will see that 
the vertex with label ``1'' takes  different positions in these different graceful copies of this same (isomorphic) tree 
isomorphic to the (3 star-1 path-tree) given in (ii) above. 

Because of the presence of label ``1'' at all inequivalent vertices (i.e. vertices belonging to different equivalent sets) 
we can generate all possible graceful trees on six vertices, corresponding to all possible non-isomorphic trees one can get 
by appending vertex-pair $(1,6),$ from (3 star-1 path-tree) given in (ii) above.

{\bf (c)} The graceful trees given by lattice paths (3),(6),(7) and (10) are all isomorphic to (4 path-tree) given in (iii) 
above. 

If we will draw corresponding tree diagrams from the above lattice paths, (3),(6),(7) and (10), we will see that 
the vertex with label ``1'' takes different positions in these different graceful copies of this same (isomorphic) tree   
isomorphic to a (4 path-tree) given in (iii) above. 

Because of the presence of label ``1'' at all inequivalent vertices (i.e. vertices belonging to different equivalent sets) 
we can generate all possible graceful trees on six vertices, corresponding to all possible non-isomorphic trees one can get 
by appending vertex-pair $(1,6),$ from isomorphic to a (4 path-tree) given in (iii) above. 

{\bf Definition 4.2:} A right-shifted $n$-lattice is the $n$-lattice that results by shifting the vertex labels in the 
vertex-pairs in the $n$-lattice by unity, i.e. we change vertex labels $i \to (i+1),$ for all $i \in \{1,2,\cdots,n\}.$ 

Thus, the so called right-shifted $n$-lattice is the one as shown below:
   \begin{align*}
    (2,3)(3,4)(4,5)(5,6)\cdots(n,n+1)     
   \end{align*}
   \begin{align*}
   (2,4)(3,5)(4,6)\cdots(n-1,n+1)
   \end{align*}
   \begin{align*}
   (2,5)(3,6)\cdots (n-2,n+1)
   \end{align*}
   \begin{align*}
   \vdots 
   \end{align*}
   \begin{align*}
   (2,n)(3,n+1)
   \end{align*}
   \begin{align*}
   (2,n+1) 
   \end{align*}
   
{\bf Observation 4.6:} By mapping each vertex label $i\to (i+1)$ in each of the vertex-pairs of a lattice path, $L,$ 
belonging to $n$-lattice we get a lattice path, $RL,$ in the right-shifted $n$-lattice. The graphs represented by $L$ 
and $RL$ are isomorphic.     

{\bf Observation 4.7:} Both, the $n$-lattice and the right-shifted $n$-lattice, exist as sub-lattices of the $(n+1)$-lattice. 

{\bf Definition 4.3:} A down-shifted $n-lattice$ is a sub-lattice of an $(n+1)$-lattice obtained by deleting the first 
row of an $(n+1)$-lattice.

{\bf Theorem 4.6:} The lattice paths in the $n$-lattice and the right-shifted $n$-lattice are exactly identical in number 
and content, i.e. for every lattice path, $L,$ belonging to $n$-lattice there exists the lattice path, $RL$ in the 
right-shifted $n$-lattice under mapping of vertex label $i\to (i+1)$ and they represent isomorphic copies of same graph. 

{\bf Proof:} Obvious. The lattice paths $L,RL$ are pictorially exactly identical. The only difference is in the labeling 
of the vertices. The only difference in the lattice path $L,$ belonging to $n$-lattice and the lattice path, $RL$ in the 
right-shifted $n$-lattice obtained under mapping of vertex label $i\to (i+1)$ is that in $L$ the vertex 
labels $\in \{1,2,\cdots,n\}$ while in $RL$ the vertex labels $\in \{2,3\cdots,(n+1)\}.$ 

{\bf Observation 4.8:} The lattice paths in the so called down-shifted $n$-lattice, a sub-lattice of an $(n+1)$-lattice 
obtained by deleting (ignoring) the first row of an $(n+1)$-lattice, contains at least two copies of $(n,n-1)$-trees of each 
isomorphism type whose edges are labeled as $k,$ $k \in \{2,3,\cdots,n\}.$     

{\bf Observation 4.9:} Every (unlabeled) tree can be viewed as (can be looked upon as) certain juxtapositions (joins) of stars 
(star subtrees) and path-trees (path subtrees). 

{\bf Examples:}

\begin{tikzpicture}[node distance={15mm}, thick, main/.style = {draw, circle}]

\node[main] (1) {$x_1$}; 
\node[main] (2) [above right of=1] {$x_2$}; 
\node[main] (3) [below right of=1] {$x_3$}; 
\node[main] (4) [above right of=3] {$x_4$}; 
\node[main] (5) [above right of=4] {$x_5$}; 
\node[main] (6) [below right of=4] {$x_6$};
\node[main] (7) [below right of=5] {$x_7$}; 
\draw (1) -- (2); 
\draw (1) -- (3);
\draw (1) -- (4);
\draw (4) -- (5);
\draw (4) -- (6);
\draw (4) -- (7);
\end{tikzpicture}

\begin{tikzpicture}[node distance={15mm}, thick, main/.style = {draw, circle}] 

\node[main] (1) {$x_1$}; 
\node[main] (2) [above right of=1] {$x_2$}; 
\node[main] (3) [below right of=1] {$x_3$}; 
\node[main] (4) [above right of=3] {$x_4$}; 
\node[main] (5) [above right of=4] {$x_5$}; 
\node[main] (6) [below right of=4] {$x_6$};
\node[main] (7) [below right of=5] {$x_7$}; 
\draw (1) -- (2); 
\draw (1) -- (3);
\draw (1) -- (4);
\draw (4) -- (5);
\draw (4) -- (6);
\draw (5) -- (7);
\end{tikzpicture}

In the first example we have a subtree which is a 4-star (a star tree with four edges) i.e. a $(5,4)$-tree made up of edges 
$(x_1,x_4),(x_4,x_5),(x_4,x_6),(x_4,x_7)$ and two 1-paths (two path subtrees containing one edge) made up of 
edges $(x_1,x_2),(x_1,x_3)$ and these trees are juxtaposed (joined) at the vertex with label $x_1.$ 

One can represent this tree in terms of stars and paths as follows: 

\begin{tikzpicture}[node distance={30mm}, thick, main/.style = {draw, circle}] 
\node[main] (1) {$4-Star$}; 
\node[main] (2) [above left of=1] {$1-Path$}; 
\node[main] (3) [below left of=1] {$1-Path$}; 
\draw (1) -- (2); 
\draw (1) -- (3);
\end{tikzpicture}

We call such a representation of a tree in terms of stars and paths as {\bf star-path representation} of that tree. It is 
not unique and one can have several star-path representations as per one's choice.

In the second example we have a subtree which is a 3-star (a star tree with three edges) i.e. a $(4,3)$-tree made up of edges 
$(x_1,x_2),(x_1,x_3),(x_1,x_4)$ and a 1-path (a path subtree containing one edge, $(4,6)$) and a 2-path 
(a path subtree containing two edges) made up of edges $(x_4,x_5),(x_5,x_7)$ and these trees are juxtaposed (joined) 
at the vertex with label $x_4.$ 

One can represent this tree in terms of a star-path representation as follows: 

\begin{tikzpicture}[node distance={30mm}, thick, main/.style = {draw, circle}] 
\node[main] (1) {$3-Star$}; 
\node[main] (2) [above right of=1] {$2-Path$}; 
\node[main] (3) [below right of=1] {$1-Path$}; 
\draw (1) -- (2); 
\draw (1) -- (3);
\end{tikzpicture} 

In the first of the above examples a 4-star is joined (juxtaposed) to two 1-paths, thus there are {\bf two joins}.

In the second of the above examples a 3-star is joined (juxtaposed) to a 2-path, and a 1-path, thus again there 
are {\bf two joins}.

{\bf Observation 4.10:} From the star-path representation of a tree, $T$ say, one can determine the number of vertices, 
$n$ say, present in $T$ by the following formula: 

Suppose a star-path representation is obtained for a tree, $T,$ say. Suppose the total number of stars and paths used 
in the star-path representation for tree $T$ are as follows:

(i) Suppose one used $i_1$ number of $j_1$-stars, $i_2$ number of $j_2$-stars, $\cdots, i_m$ number of $j_m$-stars. 

(ii) Suppose one used $k_1$ number of $l_1$-paths, $k_2$ number of $l_2$-paths, $\cdots, k_n$ number of $l_n$-paths. 

(iii) Suppose the total number of {\bf joins} one used is $p.$  

Then in this case, 

{\bf Theorem 4.7:} The total number of vertices, $n,$ present in the above tree, $T,$ are 
\begin{align*}
n = [\sum_{r = 1}^m i_r \times (j_r+1) + \sum_{s = 1}^n k_s\times (l_s+1)] - p
\end{align*}

{\bf Proof:} A $j_r$-star contains $(j_r+1)$ vertices and also an $l_s$-path contains $(l_s+1)$ vertices. If there are 
$i_r$ number of $j_r$-stars then they together contain $i_r \times (j_r+1)$ vertices for every $r \in \{1,2,\cdots,m\}$ and 
if there are $k_s$ number of $l_s$-paths then they together contain $k_s \times (l_s+1)$ vertices for 
every $s \in \{1,2,\cdots,n\}.$ Also, every join is 
nothing but a vertex where some $q$-star(s) (or path(s)) and some $t$-star(s) (or path(s)) meet therefore at each such 
join the corresponding vertex is double counted. Therefore, we have to subtract the number of joins $p$ to eliminate 
double counting of vertices. Hence etc.

{\bf Remark:} As an illustration we calculate below ``the count of vertices'' using the star-path representations 
and the formula (Theorem 4.7) for the trees given in the above examples.  

(a) The star-path representation for the tree in the first example contains one 4-star, two 1-paths, and two joins. 
Therefore, from the above formula the number of vertices in this tree are $1\times(5) + 2\times(2) - 2 = 7.$ 

(b) The star-path representation for the tree in the second example contains one 3-star, one 2-path, one 1-path, 
and two joins. Therefore, from the above formula the number of vertices in this tree 
are $1\times(4) + 1\times(3) + 1\times(2) - 2 = 7.$

{\bf Observation 4.11:} Trees thus can be represented (in a compact way) as certain juxtapositions (joins) 
of stars and paths-trees of various sizes, juxtaposed (joined) in some ways. 

{\bf Definition 4.4:} A tree, some (unlabeled) $(n,n-1)$-tree, when represented as juxtaposition (join) of some stars 
and paths in some way is called a star-path representation for that tree. It is easy to check that the star-path 
representation corresponding to a tree is not necessarily unique. In fact there can be more than one star-path representations 
for a tree of the same isomorphism type.

{\bf Definition 4.5:} A peripheral-tree is a star (or path-tree) lying on the boundary of some star-path representation of 
some (unlabeled) $(n,n-1)$-tree, and is called a peripheral-star (or a peripheral-path-tree) as a {\bf part} of that 
star-path representation.

{\bf Observation 4.12:} $n$-lattice contains, in terms of lattice paths, all possible graceful graphs as seen above. 
This $n$-lattice actually also contains (in terms of lattice paths) all possible $(n,n-1)$ graceful trees and thus 
they among them contain gracefully labeled $(n,n-1)$-trees of each isomorphism type. 

{\bf Observation 4.13:} It is possible to construct a lattice path corresponding to an (unlabeled) $(n,n-1)$-tree of 
each isomorphism type by trial and error using following steps:

(i) Find a star-path representation. Note that it is not unique and one can represent a given $(n,n-1)$-tree in 
multiple ways in terms of star-path representations. 

(ii) Out of the different stars or path-trees, present as parts, in the star-path representation of the given $(n,n-1)$-tree 
choose some peripheral-star (or peripheral-path-tree) for labeling. 

(iii) Label this chosen peripheral-star (or peripheral-path-tree) by building a partial lattice path by starting from 
the row at the bottom of an $n$-lattice by first choosing vertex pair $(1,n)$ then by appending it by a vertex-pair 
in the row above this bottom row between the vertex-pairs $\{(1,n-1),(2,n)\}.$ Continue this process of appending a 
suitable vertex-pair from next upper row to earlier chosen vertex-pairs forming a partial lattice path. If the chosen 
peripheral-tree is a $k+1$-star then thus constructed partial lattice path for this peripheral star 
will be $(1,n)(1,n-1)\cdots(1,n-k)$ or $(1,n)(2,n)\cdots(k+1,n).$ If the peripheral-tree is a $k+1$-path-tree then thus 
constructed partial lattice path for this $k+1$-path-tree will be a zigzag partial lattice path formed by choosing suitable 
vertex-pairs to append to the earlier formed partial lattice path, once from left side and then from right side, so that 
finally thus created partial lattice path will be a zigzag partial lattice path and will 
represent the chosen peripheral-tree. We can thus get a properly labeled $k+1$-star or $k+1$-path-tree, as required.  

(iv) Now continue labeling of the next star or path-tree (in the star-path representation) adjacent to the above labeled 
peripheral-star (or peripheral-path-tree) by suitably extending the above formed partial lattice path further and so on. 

(v) Complete the formation of an entire lattice path representing some graceful graph (this completes when the extension 
of partial lattice path terminates with the appending of some vertex-pair from top row).

{\bf Observation 4.14:} All those lattice paths corresponding to trees in an $n$-lattice together contain at least two 
graceful copies (a lattice path and its mirror image taken in the axis of symmetry) for every unlabeled $(n,n-1)$-tree 
(i.e. an $(n,n-1)$-tree of each isomorphism type).

{\bf Observation 4.15:} In terms of lattice paths in an $n$-lattice we can construct all possible juxtapositions 
(joins) of stars and path-trees of all possible sizes which in tern generate all $(n,n-1)$ graceful trees, at least two, for 
each isomorphism type. 

{\bf Observation 4.16:} We can generate at least one graceful copy of for every unlabeled $(n+1,n)$-tree by extending 
all lattice paths corresponding to graceful trees in $n$-lattice by appending them by an edge $(1,n+1).$ 

{\bf Example:}

In a $5$-lattice we have all possible lattice paths of length $4.$ If we will extend all the lattice paths in $5$-lattice 
by appending the vertex-pair $(1,6)$ then through this extension we get all possible lattice paths of length $5$ and 
therefore we get graceful $(6,5)$-trees such that at least one copy of every unlabeled $(6,5)$-tree will be   
present in those lattice paths. 

{\bf Explanation:} 

By attaching $(1,6)$ to {\bf first lattice path} this first lattice path 

$(1,5)(1,4)(2,4)(2,3)$ 

which is pictorially a straight path of length $4,$ (corresponding to a $4$-star) extends to  

$(1,6)(1,5)(1,4)(2,4)(2,3)$ 

which is also a {\bf straight path} and now of length $5,$ (corresponding to a $5$-star). 

PICTORIALLY it looks like: 

\begin{tikzpicture}[node distance={15mm}, thick, main/.style = {draw, circle}] 

\node[main] (1) {$(1,2)$}; 
\node[main] (2) [below right of=1] {$(1,3)$}; 
\node[main] (3) [below right of=2] {$(1,4)$}; 
\node[main] (4) [below right of=3] {$(1,5)$}; 
\node[main] (5) [below right of=4] {$(1,6)$};  
\draw (1) -- (2); 
\draw (2) -- (3);
\draw (3) -- (4);
\draw (4) -- (5);

\end{tikzpicture} 

If we will continue extending all other lattice paths of length $4$ to get lattice path of length $5$ by attaching 
vertex-pair $(1,6)$ we get all possible new lattice paths, corresponding to all possible juxtapositions of 
stars and path-trees, of length $5$ as follows:

By attaching $(1,6)$ to {\bf second lattice path} which is a (3-star- 1 path) changes to (4-star- 1 path).

By attaching $(1,6)$ to {\bf third lattice path} which is a (4-path = 2 path-2 path) changes to (3-star- 2 path). 

By attaching $(1,6)$ to {\bf fourth lattice path} which is a (1 path-3-star) changes to (two 3-stars having 1 common edge). 

By attaching $(1,6)$ to {\bf fifth lattice path} which is a (3-star- 1 path) changes to (4-star- 1 path).

By attaching $(1,6)$ to {\bf sixth lattice path} which is a (4-path = 2 path-2 path) changes to (3-star-2 path). 

By attaching $(1,6)$ to {\bf seventh lattice path} which is a (4-path) changes to (3-star having two 1 paths attached to its 
two different end vertices).

By attaching $(1,6)$ to {\bf eighth lattice path} which is a (3-star- 1 path) changes to (two 3-stars having 1 common edge).

By attaching $(1,6)$ to {\bf ninth lattice path} which is a (3-star-1 path) changes to (3-star- 2 path). 

By attaching $(1,6)$ to {\bf tenth lattice path} which is a (4-path) changes to (5 path). 

By attaching $(1,6)$ to {\bf eleventh lattice path} which is a (3 star-1-path) changes to (3-star having two 1 paths 
attached to its two different end vertices). 

and lastly by attaching $(1,6)$ to {\bf twelveth lattice path} which is a $4$-star we get the following extended 
lattice path: 

$((1,6)(1,5)(2,5)(3,5)(4,5)$ 

which is a (4-star-1 path).

PICTORIALLY it looks like:

\begin{tikzpicture}[node distance={15mm}, thick, main/.style = {draw, circle}] 

\node[main] (1) {$(4,5)$}; 
\node[main] (2) [below left of=1] {$(3,5)$}; 
\node[main] (3) [below left of=2] {$(2,5)$}; 
\node[main] (4) [below left of=3] {$(1,5)$}; 
\node[main] (5) [below right of=4] {$(1,6)$}; 
\draw (1) -- (2); 
\draw (2) -- (3);
\draw (3) -- (4);
\draw (4) -- (5); 

\end{tikzpicture} 

Thus, the extension of all lattice paths present in $5$-lattice by appending the vertex-pair $(1,6)$ leads to 
extension of a star or a path-tree whichever is present at the bottom of that lattice path in the $5$-lattice 
gets extended, i.e. if there is a $k$-star (a straight part of length $k$) or a $k$-path-tree 
(a zigzag part of length $k$) is present as a 
part of the lattice path at the bottom of that lattice path under consideration then by appending the 
vertex-pair $(1,6)$ the $k$-star becomes either a $(k+1)$-star or a (k-star -1 path) and the $k$-path-tree 
becomes either a $(k+1)$-path-tree or ((k-2)-path-tree - 3 star). It is easy to check 
that we thus get all possible lattice paths covering in them gracefully labeled $(6,5)$-trees of each 
isomorphism type.

\section{Graceful Tree Conjecture} 

In this section we will settle the Graceful Tree Conjecture in the affirmative. We show that for every 
unlabeled $(n,n-1)$-tree there exists a lattice path which represents the same (isomorphic) tree with 
desired graceful labeling. 

As stated in Observation 4, every (unlabeled) $(n,n-1)$-tree can be represented as certain 
juxtaposition (join) of stars and path-trees and we call it the star-path representation for that $(n,n-1)$-tree. 
Let us start labeling of some star or path-tree present in that star-path representation and choose the 
corresponding part of the lattice path in the $n$-lattice. For example if there exists a $k\le n$-star as a 
part of the star-path representation and suppose we start to label it first, then we can label it 
in $C^n_k = \frac{n!}{k!(n-k)!}$ ways by selecting some $k$ pairs on the left-boundary (or right-boundary) 
of the $n$-lattice. Note that the so called left-boundary of the $n$-lattice is made up of vertex-pairs: 

$(1,n)(1,n-1)(1,n-2)(1,n-3) \cdots (1,3)(1,2)$ 

and the right-boundary of the $n$-lattice is 

$(1,n)(2,n)(3,n)(4,n) \cdots (n-2,n)(n-1,n).$ 

Note that in such selection of $k$ pairs on the left-boundary vertex with label $``1''$ is present in each choice of the 
vertex-pairs. Thus, in choosing some $k$ pairs, suppose we first choose vertex-pair $(1,i_1)$ then in this choice we have 
taken two vertices, $1$ and $i_1.$ Next, suppose we choose vertex pair $(1,i_2), i_1\neq i_2,$ then in this choice we have 
added exactly one new vertex, $i_2,$ other than the earlier chosen vertices. If we continue in this way and form a straight 
part of lattice path such that we add exactly one new vertex at each selection of a vertex-pair then we form a tree which 
is a $k$-star as follows: 

$(1,i_1)(1,i_2)(1,i_3) \cdots (1,i_{(k-1)})(1,i_k),$ such that $i_1 \neq i_2 \neq i_3 \neq \cdots \neq i_k.$ Note 
that we can proceed on similar lines and also form a $k$-star by selecting some $k$ vertex-pairs on the right-boundary and 
can form (this straight) part of lattice path.  

Now, instead of first labeling a $k$-star if we will proceed with labeling a $k$-path-tree present on the  
the star-path representation then we can proceed on similar line as above and this time form {\bf a zigzag} part of 
lattice path representing a $k$-path-tree as desired, and this can be done in several ways, in fact in more than the 
ways a $k$-star can be formed as seen above. 

Now, we will proceed to label an $l$-star (or $l$-path-tree) adjacent the above labeled $k$-star (or $k$-path-tree) at a 
vertex, which is common to both parts for maintaining adjacency, and form the further part of lattice path by observing 
that (i) while appending a new vertex-pair, each time we check that this 
vertex-pair belongs to some new unused row of the $n$-lattice and the vertices in this chosen vertex-pair do not belong 
to earlier chosen vertex-pairs during the formation earlier part of the lattice path, and (ii) exactly one new vertex 
is added while appending one new vertex-pair. 

Note that we can also proceed to label an $l$-star (or $l$-path-tree) not adjacent the above labeled $k$-star 
(or $k$-path-tree) but a part of the star-path representation of the $(n,n-1)$-tree under consideration for graceful labeling. 
In such a case two new vertices are added during the selection of a vertex-pair for the first time and then exactly one new 
vertex will be added during appending a new vertex-pair as above during completing graceful labeling for 
$l$-star (or $l$-path-tree) not adjacent the above labeled $k$-star (or $k$-path-tree). 

{\bf Example:}

Suppose we are given a, $(10,9)$-tree and suppose one can represent this tree as a star-path representation as follows: 

\begin{tikzpicture}[node distance={30mm}, thick, main/.style = {draw, circle}] 
\node[main] (1) {$4-Star$}; 
\node[main] (2) [above right of=1] {$2-Path$}; 
\node[main] (3) [below right of=1] {$3-Star$}; 
\draw (1) -- (2); 
\draw (1) -- (3);
\end{tikzpicture} 

and further suppose the $2$-Path and the $3$-Star have no common vertex, i.e. they are connected to the $4$-Star at different 
vertices. It is easy to check that a lattice path is as given below for this tree which gives the graceful representation 
for this tree. 

$(1,10)(1,9)(2,9)(2,8)(2,7)(2,6)(3,6)(3,5)(3,4)$ 

which pictorially looks as below: 

\begin{tikzpicture}[node distance={15mm}, thick, main/.style = {draw, circle}] 

\node[main] (1) {$(3,4)$}; 
\node[main] (2) [below right of=1] {$(3,5)$}; 
\node[main] (3) [below right of=2] {$(3,6)$}; 
\node[main] (4) [below left of=3] {$(2,6)$}; 
\node[main] (5) [below right of=4] {$(2,7)$}; 
\node[main] (6) [below right of=5] {$(2,8)$}; 
\node[main] (7) [below right of=6] {$(2,9)$}; 
\node[main] (8) [below left of=7] {$(1,9)$}; 
\node[main] (9) [below right of=8] {$(1,10)$}; 
\draw (1) -- (2); 
\draw (2) -- (3);
\draw (3) -- (4);
\draw (4) -- (5); 
\draw (5) -- (6); 
\draw (6) -- (7);
\draw (7) -- (8);
\draw (8) -- (9);

\end{tikzpicture} 

{\bf Theorem 5.1:} The collection of all the lattice paths corresponding to trees in an $n$-lattice together contain 
graceful copies for $(n,n-1)$-trees and that of two component forests which together contain $n$ vertices and $(n-1)$ edges 
of every isomorphism type. In other words, for every (unlabeled) $(n,n-1)$-tree and for every two component forest which 
together contains $n$ vertices and $(n-1)$ edges there exists a lattice path in the $n$-lattice that represents 
a graceful copy of that $(n,n-1)$-tree and that two component forest which 
together contains $n$ vertices and $(n-1)$ edges.

{\bf Proof:} We proceed by induction on $n.$ 

{\bf Step I:} We first settle the result for $n = 1,2,3,4,5.$

For $n = 1$ the $n$-lattice is empty and nothing to prove.  

Let $n = 2.$ The $n$-lattice in this case is 

\begin{align*}
 (1,2)
\end{align*} 

and there is only one lattice path $(1,2)$ containing this only vertex-pair which represents the only graceful tree 
containing two vertices and one edge. 

Let $n = 3.$ The $n$-lattice in this case is 

\begin{align*}
    (1,2)(2,3)  
   \end{align*}
   \begin{align*}
   (1,3)
   \end{align*} 
   
and it contains two lattice paths $(1,3)(1,2)$ and $(1,3)(2,3)$ and they both represent the only tree up to isomorphism,  
containing three vertices and two edges, in two different graceful labellings. 

Let $n = 4.$ The $n$-lattice in this case is 

\begin{align*}
    (1,2)(2,3)(3,4)  
   \end{align*}
   \begin{align*}
   (1,3)(2,4)
   \end{align*}
   \begin{align*}
    (1,4)
   \end{align*}
   
and it contains four lattice paths representing graceful trees $((1,4)(1,3)(1,2),$ $(1,4)(1,3)(2,3),$ $(1,4)(2,4)(2,3)$ 
and $(1,4)(2,4)(3,4).$ The first and the fourth of the lattice paths represent gracefully labeled (3 star) while the 
second and third lattice path represent gracefully labeled (3 path-tree). A 3-star and a 3 path-tree are the only trees 
up to isomorphism, containing four vertices and three edges.  

Let $n = 5.$ The $n = 5$-lattice in this case is 

\begin{align*}
    (1,2)(2,3)(3,4)(4,5)     
   \end{align*}
   \begin{align*}
   (1,3)(2,4)(3,5)
   \end{align*}
   \begin{align*}
   (1,4)(2,5)
   \end{align*}
   \begin{align*}
   (1,5)
   \end{align*}
   
We get first six lattice paths as follows:

$(1,5)(1,4)(1,3)(1,2)$ 

$(1,5)(1,4)(1,3)(2,3)$

$(1,5)(1,4)(2,4)(2,3)$

$(1,5)(1,4)(2,4)(3,4)$

$(1,5)(1,4)(3,5)(1,2)$

$(1,5)(1,4)(3,5)(2,3)$
   
The first lattice path in this lattice is  

$(1,5)(1,4)(1,3)(1,2)$

PICTORIALLY it looks like: 

\begin{tikzpicture}[node distance={15mm}, thick, main/.style = {draw, circle}] 

\node[main] (1) {$(1,2)$}; 
\node[main] (2) [below right of=1] {$(1,3)$}; 
\node[main] (3) [below right of=2] {$(1,4)$}; 
\node[main] (4) [below right of=3] {$(1,5)$};  
\draw (1) -- (2); 
\draw (2) -- (3);
\draw (3) -- (4);

\end{tikzpicture} 

The next lattice path in this lattice is  

$(1,5)(1,4)(1,3)(2,3)$ 

PICTORIALLY it looks like:

\begin{tikzpicture}[node distance={15mm}, thick, main/.style = {draw, circle}] 

\node[main] (1) {$(2,3)$}; 
\node[main] (2) [below left of=1] {$(1,3)$}; 
\node[main] (3) [below right of=2] {$(1,4)$}; 
\node[main] (4) [below right of=3] {$(1,5)$};  
\draw (1) -- (2); 
\draw (2) -- (3);
\draw (3) -- (4);

\end{tikzpicture}

The next lattice path in this lattice is  

$(1,5)(1,4)(2,4)(2,3)$ 

PICTORIALLY it looks like: 

\begin{tikzpicture}[node distance={15mm}, thick, main/.style = {draw, circle}] 

\node[main] (1) {$(2,3)$}; 
\node[main] (2) [below right of=1] {$(2,4)$}; 
\node[main] (3) [below left of=2] {$(1,4)$}; 
\node[main] (4) [below right of=3] {$(1,5)$};  
\draw (1) -- (2); 
\draw (2) -- (3);
\draw (3) -- (4);

\end{tikzpicture}

The next lattice path in this lattice is  

$(1,5)(1,4)(2,4)(3,4)$ 

PICTORIALLY it looks like: 

\begin{tikzpicture}[node distance={15mm}, thick, main/.style = {draw, circle}] 

\node[main] (1) {$(3,4)$}; 
\node[main] (2) [below left of=1] {$(2,4)$}; 
\node[main] (3) [below left of=2] {$(1,4)$}; 
\node[main] (4) [below right of=3] {$(1,5)$};  
\draw (1) -- (2); 
\draw (2) -- (3);
\draw (3) -- (4);

\end{tikzpicture}

The next lattice path in this lattice is  

$(1,5)(1,4)(3,5)(1,2)$ 

PICTORIALLY it looks like: 

\begin{tikzpicture}[node distance={15mm}, thick, main/.style = {draw, circle}] 

\node[main] (1) {$(1,2)$}; 
\node[main] (2) [below right of=1] {$(3,5)$}; 
\node[main] (3) [below left of=2] {$(1,4)$}; 
\node[main] (4) [below right of=3] {$(1,5)$};  
\draw (1) -- (2); 
\draw (2) -- (3);
\draw (3) -- (4);

\end{tikzpicture}

The next lattice path in this lattice is  

$(1,5)(1,4)(3,5)(2,3)$ 

PICTORIALLY it looks like: 

\begin{tikzpicture}[node distance={15mm}, thick, main/.style = {draw, circle}] 

\node[main] (1) {$(2,3)$}; 
\node[main] (2) [below right of=1] {$(3,5)$}; 
\node[main] (3) [below left of=2] {$(1,4)$}; 
\node[main] (4) [below right of=3] {$(1,5)$}; 
\draw (1) -- (2); 
\draw (2) -- (3);
\draw (3) -- (4);

\end{tikzpicture} 

Continuing on these lines we now get the following six lattice paths, from seventh to twelveth lattice paths, which are the 
mirror images of the above six lattice paths (the mirror images are taken in the line of symmetry) as follows:  

$(1,5)(2,5)(1,3)(3,4)$ 

$(1,5)(2,5)(1,3)(4,5)$

$(1,5)(2,5)(2,4)(2,3)$

$(1,5)(2,5)(2,4)(3,4)$

$(1,5)(2,5)(3,5)(3,4)$

$(1,5)(2,5)(3,5)(4,5)$

PICTORIALLY the seventh lattice path $(1,5)(2,5)(1,3)(3,4)$ looks like: 

\begin{tikzpicture}[node distance={15mm}, thick, main/.style = {draw, circle}] 

\node[main] (1) {$(3,4)$}; 
\node[main] (2) [below left of=1] {$(1,3)$}; 
\node[main] (3) [below right of=2] {$(2,5)$}; 
\node[main] (4) [below left of=3] {$(1,5)$}; 
\draw (1) -- (2); 
\draw (2) -- (3);
\draw (3) -- (4); 

\end{tikzpicture}

PICTORIALLY the eighth lattice path $(1,5)(2,5)(1,3)(3,4)$ looks like: 

\begin{tikzpicture}[node distance={15mm}, thick, main/.style = {draw, circle}] 

\node[main] (1) {$(4,5)$}; 
\node[main] (2) [below left of=1] {$(1,3)$}; 
\node[main] (3) [below right of=2] {$(2,5)$}; 
\node[main] (4) [below left of=3] {$(1,5)$}; 
\draw (1) -- (2); 
\draw (2) -- (3);
\draw (3) -- (4); 

\end{tikzpicture}

PICTORIALLY the ninth lattice path $(1,5)(2,5)(1,3)(3,4)$ looks like: 

\begin{tikzpicture}[node distance={15mm}, thick, main/.style = {draw, circle}] 

\node[main] (1) {$(2,3)$}; 
\node[main] (2) [below right of=1] {$(2,4)$}; 
\node[main] (3) [below right of=2] {$(2,5)$}; 
\node[main] (4) [below left of=3] {$(1,5)$}; 
\draw (1) -- (2); 
\draw (2) -- (3);
\draw (3) -- (4); 

\end{tikzpicture}

PICTORIALLY the tenth lattice path $(1,5)(2,5)(1,3)(3,4)$ looks like: 

\begin{tikzpicture}[node distance={15mm}, thick, main/.style = {draw, circle}] 

\node[main] (1) {$(3,4)$}; 
\node[main] (2) [below left of=1] {$(2,4)$}; 
\node[main] (3) [below right of=2] {$(2,5)$}; 
\node[main] (4) [below left of=3] {$(1,5)$}; 
\draw (1) -- (2); 
\draw (2) -- (3);
\draw (3) -- (4); 

\end{tikzpicture}

PICTORIALLY the eleventh lattice path $(1,5)(2,5)(1,3)(3,4)$ looks like: 

\begin{tikzpicture}[node distance={15mm}, thick, main/.style = {draw, circle}] 

\node[main] (1) {$(3,4)$}; 
\node[main] (2) [below right of=1] {$(3,5)$}; 
\node[main] (3) [below left of=2] {$(2,5)$}; 
\node[main] (4) [below left of=3] {$(1,5)$}; 
\draw (1) -- (2); 
\draw (2) -- (3);
\draw (3) -- (4); 

\end{tikzpicture}

PICTORIALLY the twelveth lattice path $(1,5)(2,5)(1,3)(3,4)$ looks like: 

\begin{tikzpicture}[node distance={15mm}, thick, main/.style = {draw, circle}] 

\node[main] (1) {$(4,5)$}; 
\node[main] (2) [below left of=1] {$(3,5)$}; 
\node[main] (3) [below left of=2] {$(2,5)$}; 
\node[main] (4) [below left of=3] {$(1,5)$}; 
\draw (1) -- (2); 
\draw (2) -- (3);
\draw (3) -- (4); 

\end{tikzpicture} 

We have generated all the lattice paths and they together cover all possible graceful $(5,4)$-trees which  
contain a graceful copy for tree of each isomorphism type. It is easy to check that these lattice paths are those that are 
possible to form. 
There are completely straight lattice paths corresponding to star-trees 
and there are completely zigzag lattice paths corresponding to path-trees and there are all possible ``in between'' type of 
lattice paths corresponding to trees made up from all possible juxtapositions (joins) of star-trees and path-trees, 
and therefore, these lattice paths together contain a graceful copy for a tree of each isomorphism type. 

{\bf Step II:} We assume the result for $n$-lattice i.e. we assume by induction that if we take into consideration all possible 
lattice paths in the $n$-lattice contains then together they contain a graceful copy for 
$(n,n-1)$-trees of each isomorphism type and settle the result for $(n+1)$-lattice, i.e. we will show that $(n+1)$-lattice 
contains in terms of lattice paths a graceful copy for $(n+1,n)$-trees of each isomorphism type. 

We form an $(n+1)$-lattice. Both lattices, the $n$-lattice and the right-shifted $n$-lattice, exist as sub-lattices of 
this $(n+1)$-lattice. We generate all possible lattice paths corresponding to $(n,n-1)$-trees firstly in the $n$-lattice. It 
is easy to check that these lattice paths corresponding to $(n,n-1)$-trees will not contain the vertex with label $(n+1).$ 
We then extend these lattice paths by appending them with the vertex-pair $(1,n+1).$ Since $(n+1)$ was missing vertex 
therefore these lattice paths appended with the vertex-pair $(1,n+1)$ will correspond to certain graceful $(n+1,n)$-trees 
in the $(n+1)$-lattice. 

We then generate all possible lattice paths corresponding to $(n,n-1)$-trees in the right-shifted $n$-lattice. It 
is easy to check that these lattice paths corresponding to $(n,n-1)$-trees will not contain the vertex with label $1.$ 
We then extend these lattice paths by appending them with the vertex-pair $(1,n+1).$ Since $1$ was missing vertex 
therefore these lattice paths appended with the vertex-pair $(1,n+1)$ will correspond to certain graceful $(n+1,n)$-trees. 

We have not yet covered all those lattice paths that exist partially in the sub-lattices of $(n+1)$-lattice, i.e. those 
lattice paths which partially belong to $n$-lattice and partially belong to right shifted $n$-lattice. In such lattice paths 
both vertices $1$ as well as $(n+1)$ will be present. For such a lattice path, among these lattice paths, to correspond to a 
tree  when appended with vertex-pair $(1,n+1)$ 
such a lattice path (before appending with vertex-pair $(1,n+1)$) should correspond to a forest 
containing two components such that the vertex with label $1$ belongs to one component and the vertex with label $(n+1)$ 
belongs to the other component so that such a lattice path when appended with the vertex-pair $(1,n+1)$ will correspond to 
certain graceful $(n+1,n)$-tree. 

We now see that these three types of lattice paths together exhausts (covers) all possible lattice paths 
corresponding to graceful trees in $(n+1)$-lattice for $(n+1,n)$-trees of each 
isomorphism type. 

It is easy to check that one has formed lattice paths from completely straight lattice paths corresponding to star-trees 
and completely zigzag lattice paths corresponding to path-trees and all possible other ``in between'' type lattice paths 
representing all possible juxtapositions of stars and path-trees which together will constitute all possible $(n+1,n)$-trees 
which are made up of all possible juxtapositions (joins) of stars and path-trees of all possible sizes. 

{\bf Definition 5.1} Suppose $k$ is an edge of a tree joining vertices $u,v$. If $v$ is a pendant vertex then we call $u$ 
a prependant vertex.

Suppose there exists an $(n+1,n)$-tree of some isomorphism type for which there does not exist a lattice path 
representing its graceful copy among the three types of lattice paths mentioned above. Select and delete a pendant vertex 
of this tree and assign label $1$ to the corresponding prependant vertex. Now, the tree will be an $(n,n-1)$-tree which is 
graceful by induction. If we can find a lattice path for this $(n,n-1)$-tree in the $n$-lattice without altering 
the assigned label $1$ to the prependant vertex then we reattach the deleted pendant vertex and label it by label $(n+1)$ 
and we are done. If we cannot find a lattice path for this $(n,n-1)$-tree in the $n$-lattice without altering 
the assigned label $1$ to the prependant vertex then we try the same with some other pendant vertex i.e. 
by deleting some other pendant vertex and repeating the same steps. If we cannot find a lattice path for this $(n,n-1)$-tree 
in the $n$-lattice without altering the assigned label $1$ to the prependant vertex by deleting any pendant point and 
repeating the same steps then we choose some edge of given $(n+1,n)$-tree such that both the end vertices of this edge are 
not pendant vertices. We delete this edge and label the end vertices of this deleted edge by labels $1$ and $(n+1)$ 
respectively. By deleting this edge the $(n+1,n)$-tree breaks into a forest containing two component trees. Now, 
if we can find a lattice path corresponding to two component trees in $(n+1)$-lattice without altering the assigned labels 
$1$ and $(n+1)$ to the end vertices of the deleted edge then we will reattach the deleted edge which is equivalent 
to extending the lattice path by appending the vertex-pair $(1,n+1),$ and we are done. 
If we cannot find a lattice path corresponding to two component trees in $(n+1)$-lattice without altering the assigned labels 
$1$ and $(n+1)$ to the end vertices of the deleted edge then we try the same with some other edge i.e. by deleting some other 
edge whose end vertices are not pendant ones and repeating the same steps. If we cannot find a lattice path corresponding 
to two component trees in $(n+1)$-lattice without altering the assigned labels $1$ and $(n+1)$ to the end vertices of the 
other deleted edge or he same is true for any other edge. 
It leads to a contradiction to the induction hypothesis that we cannot find lattice path for certain $(n,n-1)$-trees 
and also for certain two component forests which together contain $n$ vertices and $(n-1)$ edges.

{\bf Theorem 5.2:} Every unlabeled $(n+1,n)$-tree can be obtained by extending some peripheral star or path-tree in 
the star-path representation of an $(n,n-1)$-tree of one or more isomorphism types.

{\bf Proof:} Consider any unlabeled $(n+1,n)$-tree. Obtain its star-path representation. Choose some peripheral star 
(or path-tree) in this representation and now delete some pendant point from this peripheral star (or path-tree). We then 
will be left with some $(n,n-1)$-tree, and thus the initially taken $(n+1,n)$-tree can be looked upon as an extension of 
some peripheral star or path-tree in the star-path representation of that $(n,n-1)$-tree. If delete some other pendant point 
from some other peripheral star (or path-tree) then it will result in some other $(n,n-1)$-tree nonisomorphic to earlier 
resulted $(n,n-1)$-tree and so on.

{\bf Theorem 5.3:} Every $(n,n-1)$-tree can be gracefully labeled in terms of a lattice path by starting its labeling with 
the labeling of some peripheral star (or path-tree) in its star-path representation.

{\bf Proof:} We proceed by induction on $n,$ the number of vertices. 

For $n = 1,2,3,4$ the the trees themselves are peripheral and so the result is obvious. 

Also, For $n = 5$ two out of three trees are again peripheral themselves so nothing to prove. 

The third tree is the tree with star-path representation as (3 star - 1 path-tree). Here, we can start labeling 
by choosing 3 star as a part of lattice path $(1,5)(1,4)(1,3)$ and then continue and label the remaining 1 path as $(2,3)$. 

Suppose by induction that every $(n,n-1)$-tree can be gracefully labeled by starting its labeling with the labeling of some 
peripheral star (or path-tree) in its star-path representation. 

We show that the same holds for every $(n+1,n)$-tree. 

Take some $(n+1,n)$-tree and find its star-path representation. Choose some peripheral star (or path-tree) and 
delete a pendant vertex from this peripheral tree. By induction you can start with labeling of the resulted peripheral 
star (or path-tree whatever is chosen) and find out graceful labeling for this $(n,n-1)$-tree in terms of a lattice path. 
We then reattach the deleted pendant vertex and label the resulted edge by $(1,n+1)$ and append this vertex-pair to the 
lattice path obtained for $(n,n-1)$-tree by induction and complete the construction of lattice path for the $(n+1,n)$-tree 
under consideration representing the desired graceful labeling for that tree.

{\bf Conclusion:} All trees are graceful!

 \end{document}